\documentclass[journal]{IEEEtran}
\ifCLASSINFOpdf
   \usepackage[pdftex]{graphicx}
   \graphicspath{{../figures/}}
   \DeclareGraphicsExtensions{.pdf,.jpeg,.png}
\else
\fi
\ifCLASSOPTIONcompsoc
 \usepackage[caption=false,font=normalsize,labelfont=sf,textfont=sf]{subfig}
\else
 \usepackage[caption=false,font=footnotesize]{subfig}
\fi
\usepackage{url}

\usepackage{capekCommands}
\usepackage{tucekCommands}

\begin{document}

\title{Density-Based Topology Optimization for Characteristic Modes Manipulation}

%
%

\author{Jonas~Tucek,
        Miloslav~Capek,~\IEEEmembership{Senior~Member,~IEEE,}
        Lukas~Jelinek
\thanks{Manuscript received XXX; revised XXX.
This work was supported by the Czech Science Foundation under
project No. 21-19025M and {by the Czech Technical University in Prague
under project SGS22/162/OHK3/3T/13}. The access to the computational infrastructure of the OP VVV funded project CZ.02.1.01/0.0/0.0/16\_019/0000765 ``Research Center for Informatics'' is also gratefully acknowledged.
}
\thanks{J.~Tucek, M.~Capek and L.~Jelinek are with the Czech Technical University in Prague, Prague, Czech Republic (e-mails: {jonas.tucek;~miloslav.capek;~lukas.jelinek}@fel.cvut.cz).}
\thanks{Color versions of one or more of the figures in this paper are available online at \href{http://ieeexplore.ieee.org}{\color{black}{http://ieeexplore.ieee.org}}.
}
}

\maketitle

\begin{abstract}
A density-based topology optimization framework is developed to manipulate characteristic modes of conducting surfaces. The adjoint sensitivity analysis provides an efficient computation of the material gradient utilized by the local optimizer to update the material distribution. The modal approach naturally separates geometry and feeding synthesis, demonstrating its ability to optimize modal quantities while maintaining computational efficiency through gradient-based updates. The framework's properties and performance are illustrated through several examples, including single-mode resonance control, modal Q-factor, and multi-mode optimization.
\end{abstract}


%
\IEEEpeerreviewmaketitle

\section{Introduction}

The design of compact and efficient antennas for modern wireless systems presents significant challenges~\cite{Fujimoto_Morishita_ModernSmallAntennas}. Although numerical optimization techniques are commonly employed~\cite{Garg_AnalyticalAndCompMethodsInElmag}, they often fail to provide optimal performance and physical insight into the underlying operating principles~\cite{Balanis1989}.

 Characteristic mode analysis (CMA)~\cite{Garbacz_1965_TCM},~\cite{HarringtonMautz_TheoryOfCharacteristicModesForConductingBodies} offers a powerful alternative design approach~\cite{lau2022characteristic} by enabling the independent study and design of the radiator's geometry and its excitation mechanism. The theory of CMA has become increasingly popular in recent decades for analysis and design~\cite{lau2022characteristic},~\cite{cabedo2007theory},~\cite{capek2022computational},~\cite{adams2022antenna}  of electrically small antennas~\cite{HarringtonMautz_ComputationOfCharacteristicModesForConductingBodies}. CMA has been applied to antenna shape synthesis~\cite{garbacz1982antenna} and to influence scattering by reactive loading~\cite{HarringtonMautz_ControlOfRadarScatteringByReactiveLoading}. In the past decade, CMA has found its application in the development of mobile handset antennas~\cite{li2014designHandsetAntennas}, platform antennas for unmanned aerial vehicles~\cite{chen2013ESAonUAV}, MIMO antennas~\cite{Manteuffel_Martens-CompactMultimodeMultielementAntennaForIndoorUWB}, or metasurfaces in multiport antennas~\cite{lin2018methodSuppressingHigherOrderModes}.


Topology optimization~\cite{BendsoeSigmund_TopologyOptimization}, \ie, optimization of spatial material distribution within a design domain, has been used in conjunction with genetic algorithms~\cite{EthierMcNamara_AntennaShapeSynthesisWithoutPriorSpecificationOfTheFeedpoint} to manipulate modal characteristics. Due to CMA properties, optimization without prior knowledge of the excitation can lead to potentially novel designs, and the feeding synthesis is performed after a structural optimization. A similar optimization framework has also been used to design multiple-input multiple-output (MIMO) antennas~\cite{YangAdams_SystematicShapeOptimizationOfSymmetricMIMOAntennasUsingCM} or pixel antennas~\cite{jiang2020pixel}. Nevertheless, such an optimization framework is suitable only for problems with a low number of design variables since the governed equation is solved once for each design variable~\cite{Sigmund_OnTheUselessOfNongradinetApproachesInTopoOptim}. Furthermore, this procedure often produces highly irregular structures, characterized by isolated material islands and point vertex connections~\cite{thiel2015pointContacts}. These complex geometries typically require significant postprocessing to ensure manufacturability~\cite{EthierMcNamara_AntennaShapeSynthesisWithoutPriorSpecificationOfTheFeedpoint}.

Density-based topology optimization~\cite{BendsoeSigmund_TopologyOptimization}, widely used in the finite element method (FEM), applies continuous relaxation of the design variable to utilize adjoint sensitivity analysis~\cite{tortorelli1994design}, which efficiently computes material sensitivities independently of the number of design variables. The method produces can handle even billions of design variables~\cite{Aage_etal_TopologyOptim_Nature2017}. Furthermore, the resulting structures are smooth and regular due to the utilization of filtering techniques~\cite{sigmund2007morphology}. However, such techniques are not yet commonly applied within the method of moments (MoM) framework~\cite{Harrington_FieldComputationByMoM}, primarily due to the challenge of dealing with dense matrices inherent to MoM, which complicates the implementation of efficient sensitivity analysis.

The density formulation emerged from the field of structural mechanics~\cite{bendsoe1988generating}, in which eigenvalue problems are also widely utilized to analyze vibrational modes and stability of structures~\cite{bazant2010stability}. Hence, topology optimization can be used to move characteristic frequencies away from driving frequency~\cite{zargham2016topology}, \eg{}, while designing an engine or to maximize the fundamental buckling load of a structure~\cite{neves1995buckling},~\cite{ferrari2020towards}.  

Since the density formulation has been successfully implemented within the popular MoM framework for electromagnetic applications~\cite{wang2017novel},~\cite{tucek2023TopOptMoM_minQ}, and has proven its efficiency, this article explores the potential of combining it with CMA to manipulate the modal quantities of structures without the necessity of defining an excitation, \ie{}, the goal of this article is to develop an automated procedure for manipulating any modal objective evaluated in terms of characteristic currents and characteristic numbers. The spatially varying material distribution is introduced to facilitate topology optimization, allowing continuous material distribution between the vacuum and the conductor. As long as the objectives are differentiable, the adjoint sensitivity analysis provides sensitivities to the local optimizer under the assumption of distinct characteristic numbers~\cite{lin2020state}. Furthermore, standard filtering techniques~\cite{tucek2023TopOptMoM_minQ} are utilized to regularize the design variable and to obtain a gray but near-binary solution (material is either vacuum or conductor).

The paper is structured as follows. \secref{sec:CMA} introduces CMA within the MoM modeling, including a modification to facilitate density-based topology optimization, proposes design parameterization, and investigates the influence of the material model on the characteristic solution. The topology optimization of the modal parameters is introduced in~\secref{sec:Manipulation}. \secref{sec:flowchart} summarizes the workflow. Section~\ref{sec:examples} demonstrates optimizations for single and multi-mode resonance, including additional constraints, highlighting the proposed technique's properties, performance, and computational complexity. The paper is concluded in~\secref{sec:conclusion}.

\section{Characteristic Modes in Method of Moments}\label{sec:CMA}
In this paper, MoM is applied to numerically model the physical behavior of radiating structures under the assumption of a time-harmonic steady state at angular frequency~$\omega$ with convention~$\partial/\partial t \to \J\omega$, with $\J$ being the imaginary unit. The associated continuous impedance operator is described by the electrical field integral equation (EFIE)~\cite{Gibson_MoMinElectromagnetics} for perfectly conducting bodies and is converted to the impedance matrix
\begin{equation}
    \M{Z} = \M{R}_0 +\J \M{X}_0.
    \label{eq:impedance matrix}
\end{equation}
Without loss of generality, we assume spatial discretization into~$T$ triangular patches and employment of $N$~Rao-Wilton-Glisson basis functions~\cite{RaoWiltonGlisson_ElectromagneticScatteringBySurfacesOfArbitraryShape}. Hence, matrices~$\M{R}_0$ and~$\M{X}_0$ in~\eqref{eq:impedance matrix} are real, symmetric, and of size~$N\times N$.

Characteristic modes are numerically obtained through generalized eigenvalue problem~\cite{HarringtonMautz_ComputationOfCharacteristicModesForConductingBodies} 
\begin{equation}
    \M{X}_0\M{I}_n = \lambda_n \M{R}_0\M{I}_n,
    \label{eq:CMA}
\end{equation}
where each characteristic current~$\M{I}_n$ and associated characteristic number~$\lambda_n$ are both real~\cite{HarringtonMautz_ComputationOfCharacteristicModesForConductingBodies}. The magnitude of characteristic number~$\lambda_n$ measures how well the associated mode radiates, \ie{}, provides the ratio of the net reactive power to the radiated power, and reveals if the mode is capacitive ($\lambda_n<0$), inductive ($\lambda_n>0$), or resonant ($\lambda_n=0$). 

Antenna design takes advantage of the physical insights provided by CMA to improve antenna performance to meet a design goal without defining an excitation~\cite{VogelEtAl_CManalysis_PuttingPhysicsBackIntoSimulation}. By optimizing the topology of an antenna, we can adjust the characteristic number to be close to zero at a desired frequency, efficiently making it resonant~\cite{cabedo2007theory}, or modify other modal parameters such as modal Q-factor~\cite{EthierMcNamara_AntennaShapeSynthesisWithoutPriorSpecificationOfTheFeedpoint} or modal radiation efficiency~\cite{Li+etal2018a}. 
The orthogonality of CMA can also be conveniently utilized to design the MIMO systems by making more modes resonant at the desired frequency~\cite{li2014designHandsetAntennas},~\cite{dong2020design}. 

Within the surface formulation of EFIE for conducting structures, the conducting material is conveniently described by the surface resistivity~$\Rs$ and the physical response of the continuous material distribution is captured in the real-valued symmetric material matrix~$\M{R}_\rho$. To apply topology optimization, the material matrix is added to the real part of the impedance matrix in~$\eqref{eq:impedance matrix}$ and introduced to the left-hand side of~\eqref{eq:CMA} as
\begin{equation}
    \left(\M{X}_0 - \J \M{R}_\rho\right)\M{I}_n = \xi_n \M{R}_0 \M{I}_n.
    \label{eq:lossy CMA}
\end{equation}
The matrix on the left-hand side is complex, and consequently, characteristic numbers~$\xi_n$ and characteristic currents~$\M{I}_n$ are also complex. Thus, the appropriate orthogonality relationships are
\begin{align}
    \M{I}_m^\herm \M{R}_0\M{I}_n &= \delta_{mn},\\
    \M{I}_m^\herm \left(\M{X}_0 - \J \M{R}_\rho \right)\M{I}_n &= \xi_n\delta_{mn},
\end{align}
where $\delta_{mn}$ is the Kronecker delta. Because of these, the orthogonality of the radiation fields is preserved\footnote{Introduction of the material matrix~$\Zrho$ to the right-hand side of~\eqref{eq:CMA} breaks the orthogonality of modal radiation patterns~\cite{HarringtonMautzChang_CharacteristicModesForDielectricAndMagneticBodies}.}. Furthermore, the characteristic numbers~$\xi_n$ have the following physical interpretation 
\begin{equation}
    \xi_n = \frac{\M{I}_n^\herm\M{X}_0\M{I}_n}{\M{I}_n^\herm\M{R}_0\M{I}_n} - \J\frac{\M{I}_n^\herm\M{R}_\rho\M{I}_n}{\M{I}_n^\herm\M{R}_0\M{I}_n} = \lambda_n - \J \delta_n,
    \label{eq:Xi}
\end{equation}
where $\delta_n$ is the modal dissipation factor. Thus, the real part of the characteristic number~$\xi_n$ has the same interpretation as in~\eqref{eq:CMA}, and the imaginary part relates the modal radiation efficiency~\cite{yla2019generalized}. 

The loss term is introduced to~\eqref{eq:CMA} as an auxiliary variable, utilized exclusively during the optimization phase to facilitate the density formulation of topology optimization. Although this may suggest a connection to real material properties, in our approach, the loss term serves as an instrument for topology modification rather than simulating real metals. By employing this auxiliary variable, we gain the flexibility to modify and control key parameters of the topology during the optimization phase. 

\subsection{Design Parametrization}
This contribution aims to distribute a perfect electric conductor (PEC) in a fixed design domain to meet the requirements on the modal quantities associated with the shape of an antenna. To facilitate the density formulation of topology optimization, each pixel of the discretization grid is associated with the continuous design variable~$\des\in[0,1]^{T\times 1}$, which is piecewise constant over the triangular grid, \ie, subscript~$_t$ represents $t$-th triangle. 

The design variable represents the density of a material, \ie{}, it is permitted to take intermediate values between the air ($\des=0$) and PEC ($\des=1$). The interpolation function from~\cite{tucek2023TopOptMoM_minQ} is utilized to project the design variable onto the surface resistivity~$\Rs(\des)$ as
\begin{equation}
    \Rs(\des) = \Zair \left(\dfrac{\Zmet}{\Zair} \right)^{\des/(2-\des)}, 
    \label{eq:Interpolation}
\end{equation}
where $\Zair$ and $\Zmet$ represent the surface resistivity of the vacuum and PEC, respectively. The physical response of the material distribution~$\Rs(\des)$ is captured in the material matrix~$\M{R}_\rho(\des)$ as
\begin{equation}
    \Zrho(\des) = \displaystyle\sum_{t=1}^T \Rs(\des)\Zrhoe,
\end{equation}
where $\Zrhoe$ is the material element matrix~\cite{tucek2023TopOptMoM_minQ}. The behavior of CMA~\eqref{eq:lossy CMA} is directly influenced by changes in~$\M{R}_\rho(\des)$ as the design variable~$\des$ is perturbed. The interplay between the design variable and the auxiliary loss term allows the optimization process to fine-tune the material distribution, enabling a convergence toward the optimized solution. 

The surface resistivity of the air and PEC are set based on numerical investigations performed by the authors, see~\appref{app:Interpolation behavior Zs}, as
\begin{align}
    \des=0 &\to \Zair = 10^5\;\Omega, \\
    \des=1 &\to \Zmet=0.01\;\Omega,
    \label{eq:Resistivities}
\end{align}
which sufficiently models the physical response of the density field during the optimization. It is important to note that these resistivity values are used solely to facilitate the optimization procedure. The goal is to obtain a near-binary solution that maintains performance when postprocessed to consist of only vacuum and PEC, which is then analyzed in lossless settings~\eqref{eq:CMA}.

\begin{figure}[!t]
\centering
\includegraphics[width=3.25in,clip]{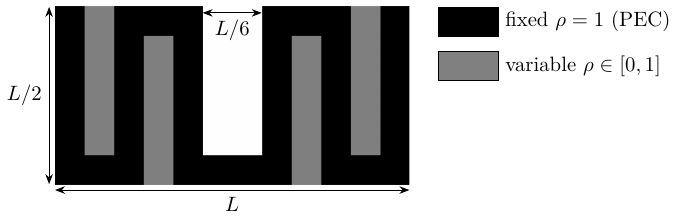}
\caption{The structure used to investigate the behavior of the material model. The region is divided into three parts: a passive region (white) with~$\des = 0$ (vacuum), a fixed region (black) with~$\des = 1$ (PEC), a design region (gray) with variable~$\rho \in[0,1]$, which translates to~$\Rs(\rho)$ through~\eqref{eq:Interpolation}. Hence, $\rho=0$ leads to a meander line antenna, while $\rho=1$ results in a slot-loaded plate.}
\label{fig:structure for material model investigations}
\end{figure}

\subsection{The Influence of the Material Model}\label{sec:Influence}
\begin{figure}[!t]
\centering
\includegraphics[width=3.25in,clip]{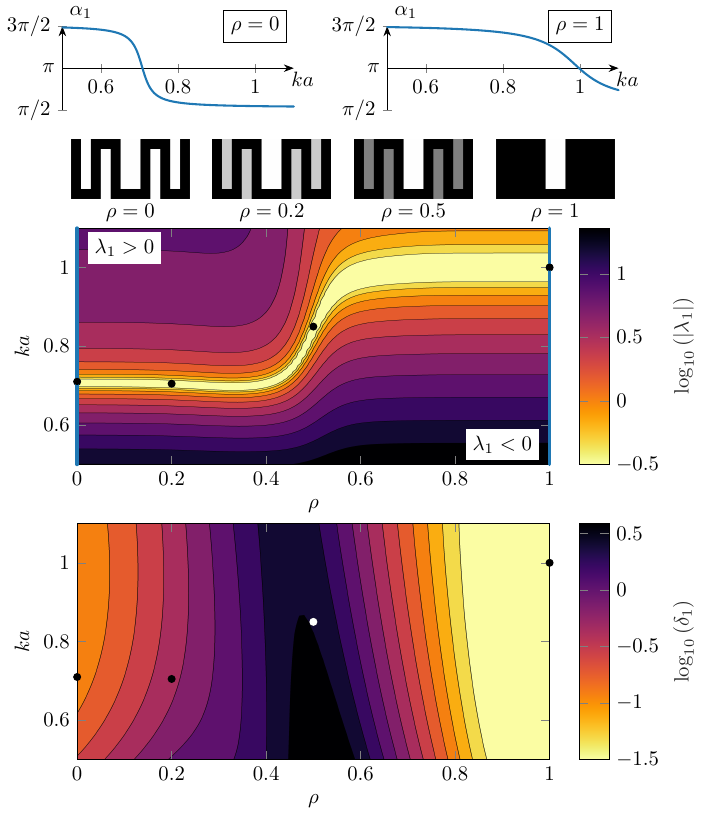}

\caption{Contour plots showing magnitudes of the both parts $\lambda_1$ and $\delta_1$ of the first characteristic number $\xi_1=\lambda_1-\J\delta_1$ as a function of electrical size~$ka$ and design variable~$\rho$ which relates the surface resistivity~$\Rs(\des)$ through~\eqref{eq:Interpolation}. The structure transitions between two geometries as $\rho$ varies from 0 to 1 and is represented by each marker, see snapshots above the axis. The minimum values of the color bars are truncated to improve the readability. The top panes show the characteristic angles~$\alpha_1$ of both extreme cases with~$\rho = \{0,1\}$ to highlight the resonance. 
}
\label{fig:material model contourplots}
\end{figure}

The influence of the material model in CMA~\eqref{eq:lossy CMA} is further examined on a continuous transition between two topologies. At the same time, the density field~$\rho$ is varied from~0 to~1 in the design region of the slot-loaded plate, see~\figref{fig:structure for material model investigations}. The first significant characteristic number~$\xi_1=\lambda_1-\J\delta_1$ (the lowest magnitude) of the structure is evaluated for each pair of electrical size~$ka$, where~$k$ is the wave number and~$a$ is the smallest circumscribing sphere, and~$\rho$ in the intervals~$ka\in[0.5,1.1]$ and~$\rho\in[0,1]$, respectively, with impedance and material matrices being evaluated in AToM~\cite{atom}. Both parts of the first characteristic number~$\lambda_1$ and $\delta_1$ are depicted in~\figref{fig:material model contourplots}.

The boundary cases in terms of~$\rho$ are discussed first. Consider the meander line antenna, \ie,~$\rho=0$, the first significant mode has the sharp resonance at $ka\approx0.7$ which is illustrated by the characteristic angle~\cite{Newman_SmallAntennaLocationSynthesisUsingCharacteristicModes}
\begin{equation}
    \alpha_n=\pi-\arctan(\lambda_n),
\end{equation}
which shows that the mode is capacitive below the resonance and inductive above the resonance, see the plot above the structure with~$\rho=0$ in~\figref{fig:material model contourplots}. If the meanders are galvanically short-circuited, we obtain the slot-loaded plate, \ie{}, $\rho=1$, with the broader resonance at $ka\approx 1$. Furthermore, the modal radiation efficiency~\cite{luo2022efficiency} is negligibly influenced in these cases above since there are virtually no losses.

The effect of the continuous interpolation of~$\Rs(\rho)$ through~\eqref{eq:Interpolation} in the design part in~\figref{fig:structure for material model investigations} is also investigated, see~\figref{fig:material model contourplots}. The real part of~$\xi_1$ (middle pane) evolves smoothly between the boundary cases. The magnitude of the imaginary part of the characteristic number~$\xi_1$ (bottom pane), \ie{}, modal dissipation delta~$\delta_1$, is extremely influenced by the extra ohmic losses induced by the intermediate values of~$\rho$.  

Regarding the improvement in resonance characteristics, the meander line antenna on the left ($\rho = 0$) exhibits a sharper resonance at a lower~$ka$ compared to the slot-loaded plate on the right ($\rho = 1$). By varying $\rho$ continuously, the optimizer can shift the resonance frequency or broaden the resonance based on the design objectives. 

When the evaluation of modal dissipation factor~$\delta_n$ is included as part of the optimization objective,~$\delta_n$ may dominate its value and hide the resonance properties, particularly in the region of intermediate~$\rho$. Thus, the optimizer must navigate through the peak at $\rho \approx 0.5$, which introduces a strong penalty in the optimization process. This penalty is referred to as self-penalization~\cite{HassanWadbroBerggren_TopologyOptimizationOfMetallicAntennas}. The optimizer can be strongly biased towards lossless designs ($\rho\to1 \ \T{or} \ 0$), and it can neglect the intermediate region of~$\rho$ that could lead to better overall performance. This creates a self-imposed penalty, as the optimizer cannot explore or settle in these intermediate regions. Nevertheless, not including~$\delta_n$ in the optimization also leads to a suboptimal design. Hence, a compromise must be found, as demonstrated in the examples.

In summary, this investigation highlights the critical role of material interpolation in influencing the characteristic spectrum when an arbitrary number of design variables are introduced. Since the general optimization task involving CMA~\eqref{eq:lossy CMA} is self-penalized, the objectives must be carefully constructed to reach the best possible solution. 

\section{Manipulation with Modal Characteristics}\label{sec:Manipulation}
Many radiating properties, \eg{}, radiation efficiency, bandwidth, or resonance, can be related to characteristic numbers. Hence, CMA can be conveniently combined with a topology optimizer without prior knowledge of feeding location~\cite{EthierMcNamara_AntennaShapeSynthesisWithoutPriorSpecificationOfTheFeedpoint}. Within our method, the general optimization task is defined as the nested formulation~\cite{BendsoeSigmund_TopologyOptimization} and reads
\begin{equation}
    \begin{aligned}
        \underset{\des}{\T{minimize}} &\quad f(\xi_n, \M{I}_n, \des)  \\
         \T{subject~to}
          &\quad h_i(\xi_n, \M{I}_n, \des) \leq 0, \\
         &\quad \left(\M{X}_0 - \J \M{R}_\rho(\des)\right)\M{I}_n = \xi_n \M{R}_0 \M{I}_n,\\
         &\quad 0\leq \des \leq 1,\quad t=\left\{1,\dotsc,T\right\},\\
         \label{eq:general TopOpt task}
    \end{aligned}
\end{equation}
where $f$ and $h_i$ are an arbitrary objective function and constraints, respectively, expressed by the design variable~$\des$, characteristic numbers~$\xi_n$ and characteristic currents~$\M{I}_n$. 

Tackling~\eqref{eq:general TopOpt task} as it stands often results in mesh-dependent designs with point vertex connections~\cite{thiel2015pointContacts}, which cannot be physically captured within the MoM models and needs additional postprocessing. Thus, the density design variable~$\des$ is filtered (smoothened) by the density filter~$\densityFilter(\cdot)$~\cite{tucek2023TopOptMoM_minQ} which is utilized to impose a length-scale~\cite{bruns2001topology} and to regularize the solution space. However, it introduces extra intermediate densities in the design, inducing extra ohmic losses, and thus, it extremely influences modal dissipation factor~$\delta_n$. Therefore, the smooth projection filter~$\heavisideFilter(\cdot)$~\cite{wang2011projection} is implemented in the $\beta$-continuation scheme~\cite{tucek2023TopOptMoM_minQ} to limit the occurrence of intermediate densities, \ie{}, filter sharpness~$\beta$ is gradually increased during the optimization to obtain a near-binary solution. The filters are applied as
\begin{equation}
    \desF = \heavisideFilter({\densityFilter({\des}})),
\end{equation}
where $\desF$ is the projected and filtered design variable. The effect of both filters is also summarized and demonstrated in~\cite{tucek2023TopOptMoM_minQ}. Including the filtering into~\eqref{eq:general TopOpt task} leads to
\begin{equation}
    \begin{aligned}
        \underset{\des}{\T{minimize}} &\quad f(\xi_n, \M{I}_n, \desF)  \\
         \T{subject~to}
          &\quad h_i(\xi_n, \M{I}_n, \desF) \leq 0, \\
         &\quad \left(\M{X}_0 - \J \M{R}_\rho(\desF)\right)\M{I}_n = \xi_n \M{R}_0 \M{I}_n,\\
         &\quad \desF = \heavisideFilter({\densityFilter({\des}})),\\
         &\quad 0\leq \des \leq 1,\quad t=\left\{1,\dotsc,T\right\},\\
         \label{eq:filters TopOpt task}
    \end{aligned}
\end{equation}
where the objectives are computed for the physically filtered design variable~$\desF$. The topology optimization task~\eqref{eq:filters TopOpt task} is solved by the gradient-based optimizer, the method of moving asymptotes (MMA)~\cite{svanberg1987method}. Considering no constraints~$h_i$ in~\eqref{eq:filters TopOpt task} the slightly modified version of MMA is utilized~\cite{da2021local},~\cite{giraldo2021polystress}. The desired number of characteristic modes~\eqref{eq:lossy CMA} is evaluated by implicitly restarted Arnoldi method~\cite{Arnoldi_ThePrincipleOfMinimizedIterationsInTheSolution}.

The sensitivities (material gradient) of the objectives and constraints with respect to the filtered design variable~$\desF$ are analytically computed by the adjoint sensitivity analysis~\cite{tortorelli1994design}, which is tailored for the properties of CMA~\cite{he2023eigenderivation}, see the complete derivation in~\appref{app:adjoint}. 

The adjoint sensitivity analysis presented in~\appref{app:adjoint} assumes distinct characteristic numbers. This poses potential challenges when characteristic numbers cross in the complex space parameterized by the frequency. This crossing results in the coefficient matrix in~\eqref{eq:Adjoint equation} and~\eqref{eq:adjoint lambda} becoming singular. Consequently, the continuous crossings~\cite{SchabBernhard_GroupTheoryForCMA} occurring in structures with non-Abelian point-group symmetries~\cite{McWeeny_GroupTheory} are unwanted and need to be remedied by cautiously defining the design region, \eg{}, the discretization can be non-uniform~\cite{2021_Capek_TAP_SymmetriesAndBounds}. While point crossing~\cite{Maseketal_ModalTrackingBasedOnGroupTheory} is theoretically possible, it's highly unlikely during optimization due to numerical factors causing non-symmetric distribution of intermediate densities in the design domain. Furthermore, since CMA~\eqref{eq:lossy CMA} produces complex characteristic numbers, crossings would need to occur simultaneously in both real and imaginary parts.

\begin{figure}[!t]
\centering
\includegraphics[width=3.25in,clip]{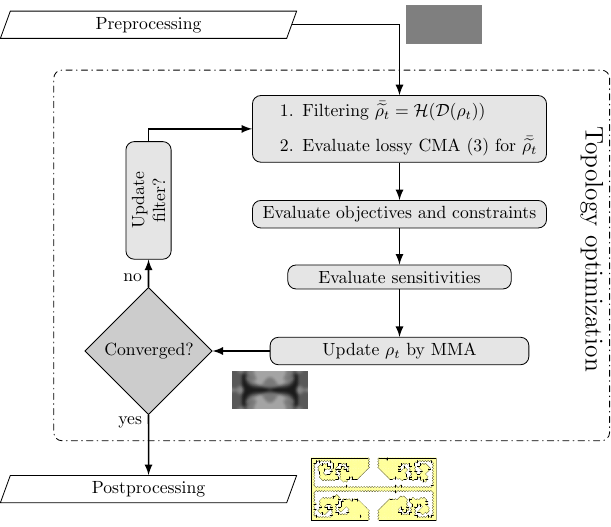}
\caption{Flowchart adapted from~\cite{tucek2023TopOptMoM_minQ} illustrates the density-based topology optimization process for manipulating characteristic modes.}
\label{fig:flowchart}
\end{figure}

\section{Optimization Workflow}\label{sec:flowchart}

The detailed workflow for density-based topology optimization manipulating with characteristic modes is presented in~\figref{fig:flowchart}. The flowchart is adapted from~\cite[Fig. 5]{tucek2023TopOptMoM_minQ} and is overall similar beside one step. It is divided into three phases: preprocessing, optimization, and postprocessing. 

The preprocessing phase involves structure discretization, computation of MoM matrices, and setup of optimization parameters such as initial density distribution, see~\cite[Section~IV]{tucek2023TopOptMoM_minQ}. The setup of the $\beta$-continuation scheme remains the same as in~\cite{tucek2023TopOptMoM_minQ}.  

The topology optimization stage begins with density and projection filtering of the design variable and with the evaluation of CMA~\eqref{eq:lossy CMA} for a desired number of modes (necessary constraints are imposed, see~\appref{app:adjoint}). Next, the objectives and constraints are computed based on characteristic numbers and currents from the previous step. The adjoint sensitivity analysis calculates the sensitivities for each objective and constraint. MMA updates the design variable~$\des$. If the convergence criteria are not met, the optimization continues with the check if the filter should be updated, and the cycle repeats, possibly with the new filter settings. The convergence and filter update criteria are the same as in~\cite[Section~IV]{tucek2023TopOptMoM_minQ}, but the projection filter is updated every 75 iterations if not stated otherwise. This process is repeated iteratively until the optimization converges to the local minimum. 

The last stage is postprocessing, \ie{}, a near-binary structure with intermediate values of~$\des$ (gray-colored design) is converted to a binary result ($\des\in\{0,1\}^{T\times 1}$) by hard thresholding~\cite{tucek2023TopOptMoM_minQ}, visualized (yellow-colored design in~\figref{fig:flowchart}) and its performance is evaluated in lossless settings~\eqref{eq:CMA}. The optimization settings, material model properties, or the choice of objectives significantly impact the result of the optimization procedure. Hence, it is advised to perform multiple runs with various initial settings~\cite{tucek2023TopOptMoM_minQ}. This is possible due to the superior scalability of the method, as demonstrated in~\secref{sec:resonance}, if compared to the optimization with genetic algorithms.

\section{Examples}\label{sec:examples}

To demonstrate the effectiveness of the proposed method, we present several benchmark examples with common objectives encountered in CMA community, \ie, the goal is to enhance the modal characteristics of a structure by optimizing the density distribution of a conductive material. The optimization procedure is implemented in MATLAB~\cite{matlab} according to~\secref{sec:flowchart}. The MoM matrices are evaluated in AToM~\cite{atom}. The examples are computed on a computer with an Intel Xeon Gold 6244 CPU (16 physical cores, 3.6\,GHz) with 384\,GB RAM. The extensive parametric studies were performed on RCI cluster~\cite{rci} on compute nodes with AMD EPYC 7543 CPU (64 physical cores, 3.1\,GHz) and with 1\,TB RAM.

The topology optimization is initiated with uniform density distribution~$\des=0.5$, and the design is subject to the density filter~$\densityFilter(.)$ with the fixed radius~$\Rmin=0.1a$ if not otherwise stated.
 
\subsection{Resonance Tuning}\label{sec:resonance}
The common objective is to enforce characteristic mode to be resonant~\cite{SafinManteuffel_ManipulationOfCWMbyImpedanceLoading}. Within the developed framework, this task can be formulated as
\begin{equation}
    \begin{aligned}
        \underset{\des}{\T{minimize}} &\quad |\lambda_1|^2 + \nu |\delta_1|^2  \\
         \T{subject~to} 
         &\quad \left(\M{X}_0 - \J \M{R}_\rho(\desF)\right)\M{I}_1 = \xi_1 \M{R}_0 \M{I}_1,\\
          &\quad \desF = \heavisideFilter({\densityFilter({\des}})),\\
         &\quad 0\leq \des \leq 1,\quad t=1,\dotsc,T,\\
         \label{eq:MakeModeResonant1}
    \end{aligned}
\end{equation}
where we manipulate the first (most significant) characteristic number, and where $\nu$ is a heuristic parameter used to change the significance of losses in the objective function. Based on our numerical investigations, we recommend starting with low values of $\nu$ and gradually increasing towards $\nu=1$, conducting multiple optimization runs to find the best solution. We advise against using $\nu$ values greater than 1, as this would induce extra ohmic losses and, therefore, excessive amplification of the self-penalization. The optimization task~\eqref{eq:MakeModeResonant1} is solved for local optimality by MMA at electrical size~$ka=0.5$, on the rectangular design region of size $L \times L/2$ which is uniformly discretized into 5000 triangles (7425 basis functions). 


\begin{figure}[!t]
\centering
\includegraphics[width=3.25in]{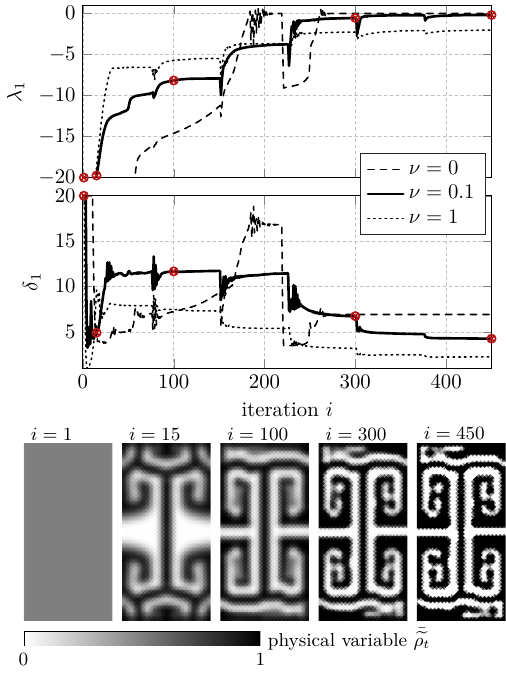}
\caption{The convergence history of three runs of the topology optimization task~\eqref{eq:MakeModeResonant1} with different values of heuristic parameter~$\nu$ included in the objective. The parameter~$\nu$ influences the magnitude of the imaginary part of the characteristic number and allows the optimizer to converge to different results. A few snapshots of the density distribution from selected iterations of the run with~$\nu=0.1$ are included below the convergence plots and highlighted by the red markers.}
\label{fig:Convergence_MakeModeResonant}
\end{figure}

Three optimization runs were performed with different values of~$\nu=\{0,0.1,1\}$, with each iteration taking approximately 40 seconds. The convergence history in~\figref{fig:Convergence_MakeModeResonant} shows both parts of the characteristic number~$\xi_1$ to demonstrate how~$\nu$ influences the optimization:
\begin{itemize}
    \item $\nu=1$ shows steep convergence but fails to achieve resonance ($\lambda_1\neq 0$) because the optimizer is driven by the self-penalization property to a local minimum with low-value of~$\delta_1$.
    \item $\nu=0$ shows the slowest and the most unstable convergence, achieving the resonance but with no control over~$\delta_1$.      
    \item $\nu=0.1$ shows gradual convergence and reaches the resonance at the expense of a slightly higher value of~$\delta_1$. Snapshots of the density distribution in several iterations are also included in~\figref{fig:Convergence_MakeModeResonant}.
\end{itemize}

The proposed method is an inherently local procedure that provides relatively smooth gray solutions near resonance, potentially eliminating point vertex connections~\cite{thiel2015pointContacts}. Smoothness is achieved through density filtering, which regularizes the solution space. This filtering process ensures solution smoothness and helps control the optimized structures' geometric complexity.
  
\begin{figure}[!t]
\centering
\includegraphics[width=3.25in]{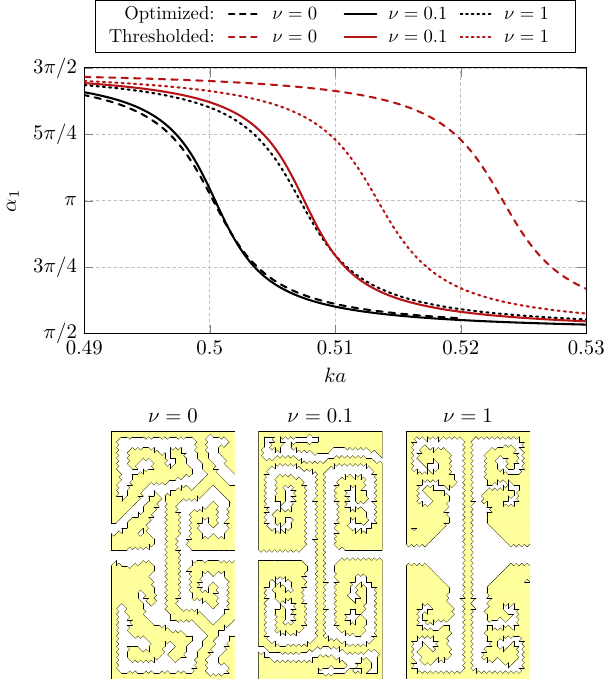}
\caption{The frequency sweep of of characteristic angle~$\alpha_1$ for all three designs optimized at electrical size~$ka=0.5$ (black lines). These near-binary solutions are hard-thresholded, impacting performance as shown by the red lines. The thresholded designs are shown as the bottom panes.}
\label{fig:Threshold_MakeModeResonant}
\end{figure}

The continuation scheme of the projection filter is utilized during the optimization, achieving a near-binary solution. The update of the projection filter results in the sudden jumps in the convergence, see~\figref{fig:Convergence_MakeModeResonant}. The near-binary solutions are converted to binary designs through hard-thresholding~\cite{tucek2023TopOptMoM_minQ}. Thresholding process inevitably causes a shift in the resonant frequency~\cite{diaz2010topology}, as shown in~\figref{fig:Threshold_MakeModeResonant}, where the characteristic angle~$\alpha_1$ is evaluated for both optimized and binary designs. It is important to emphasize that the optimized designs possess a non-negligible value of~$\delta_1$, but the thresholded binary designs are lossless ($\delta_1=0$). The binary design obtained with~$\nu=0.1$ exhibits the slightest frequency shift from the target~$ka$, suggesting that moderating the influence of losses during optimization through~$\nu$ leads to binary solutions closer to the resonance, see the structures in~\figref{fig:Threshold_MakeModeResonant}. This result indicates that properly balancing the influence of losses through~$\nu$ is crucial for maintaining performance after thresholding. The framework offers flexibility by allowing various optimization schemes, \eg{}, a continuation scheme in~$\nu$, or initializing design variables based on previous runs to reach a satisfactory solution.
 
\begin{figure}[!t]
\centering
\includegraphics[width=3.25in]{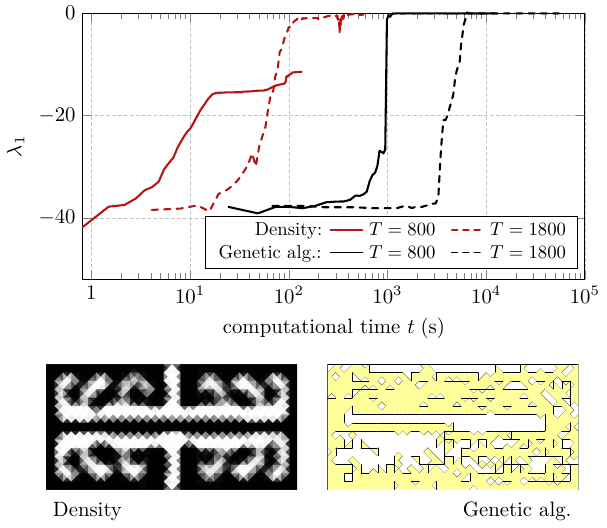}
\caption{Computational time of two topology optimization schemes (genetic algorithm~\cite{Ethier_AntennaShapeSyntesisUsingCMconcepts} and the density formulation) necessary to obtain a solution with the first significant characteristic mode close to the resonance. Results are shown for two discretization levels with $T$ triangles. Optimized designs at the $T=1800$ grid are included below for both strategies.}
\label{fig:ComputationalTime}
\end{figure}

The computational efficiency of the proposed method is evaluated against binary topology optimization performed by the genetic algorithm~\cite{Ethier_AntennaShapeSyntesisUsingCMconcepts}, see~\figref{fig:ComputationalTime}. The comparison focuses on optimizing the first characteristic mode resonance at~$ka$ using two discretization levels (800 and 1800 triangles). The density filter radius is reduced to~$\Rmin=0.05a$ to accommodate coarse discretization. The density-based optimization, implemented with~$\nu=0.1$, shows varying performance across discretization levels. While the coarser grid struggles due to ineffective density filtering (red solid line), the finer discretization successfully locates near-resonant gray solutions through proper solution space regularization~\cite{tucek2023TopOptMoM_minQ} (red dashed line).

The genetic algorithm~\cite{Ethier_AntennaShapeSyntesisUsingCMconcepts}, utilizing 200 agents, successfully finds resonant binary solutions but exhibits slower convergence compared to the density-based approach. Despite the utilization of regularization metrics~\cite{capek2021topoSensWGeom} in the objective during optimization, the genetic algorithm often produces highly irregular structures, making manufacturing challenging. Additionally, genetic algorithms struggle to effectively handle constraints in the optimization process. Consequently, both methods require postprocessing: thresholding for density-based designs and refinement for the irregular genetic algorithm solutions~\cite{Ethier_AntennaShapeSyntesisUsingCMconcepts}, see both designs in~\figref{fig:ComputationalTime}. The density-based method's superior scalability, gradient-driven convergence, and ability to handle constraints make it particularly suitable for repetitive optimization across various settings~\cite{tucek2023TopOptMoM_minQ}.

It's important to note that optimizing for resonance alone can yield numerous possible solutions. This diversity of structures offers opportunities for further optimization under additional design criteria, as demonstrated in the next section.

\begin{figure}[!t]
\centering
\includegraphics[width=3.25in]{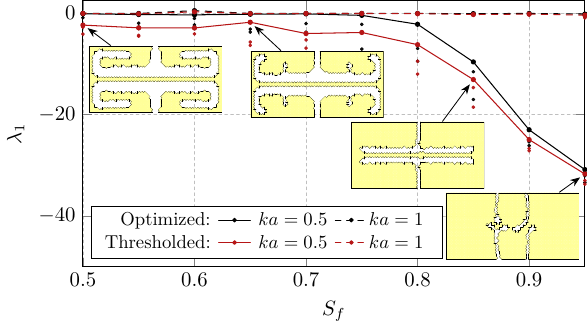}
\caption{The trade-off between the real part of the first characteristic number~$\lambda_1$ and minimum area fraction~$\Sf$ (solid and dashed lines), constructed from multiple optimization runs (individual points). The results are compared at two electrical sizes~$ka = 0.5$ and~$ka = 1$, where black curves represent optimized gray designs and red curves show their thresholded binary counterparts. Thresholding the optimized solution results in a slight drop in the objective, but the design remains close to the resonance. Thresholded binary designs are also displayed for selected area fractions, demonstrating structural evolution as~$\Sf$ increases.}
\label{fig:Cost of resonance}
\end{figure}

\subsection{Resonance Tuning Under Area Constraint}
This section introduces an additional constraint to select specific solutions from the many potential resonant structures. Surface area is a common design constraint in antenna engineering, and thus, a relationship between achieving resonance and the area covered by the material is investigated, exploring potential trade-offs between the first significant characteristic number and the surface area of a solution. This constraint on the amount of material in the design region is incorporated into the optimization task as


\begin{equation}
    \begin{aligned}
        \underset{\des}{\T{minimize}} &\quad |\lambda_1|^2 + \nu |\delta_1|^2  \\
         \T{subject~to} 
         &\quad \left(\M{X}_0 - \J \M{R}_\rho(\desF)\right)\M{I}_1 = \xi_1 \M{R}_0 \M{I}_1,\\
         &\quad \dfrac{1}{A_0}\sum_{t=1}^T\desF A_t \geq \Sf,\\
          &\quad \desF = \heavisideFilter({\densityFilter({\des}})),\\
         &\quad 0\leq \des \leq 1,\quad t=1,\dotsc,T,\\
         \label{eq:ReverseAreaConstraint}
    \end{aligned}
\end{equation}
where~$A_0$ is the total design region area,~$A_t$ is the area of the associated~$t$-th triangle, and~$\Sf$ is the minimum fraction of the design domain area. This extra constraint is parametrized by the purely heuristic parameter~$\Sf$, which restricts the minimum amount of the material. The goal, the design region, and optimization settings remain the same as in the previous section. The optimization task~\eqref{eq:ReverseAreaConstraint} is solved for several values of relative area~$\Sf$ at electrical sizes $ka=\{0.5,1\}$ using the RCI cluster~\cite{rci}.

The results are shown in~\figref{fig:Cost of resonance}, where the black and red markers represent the real part of the first characteristic number~$\lambda_1$ for optimized (gray) and thresholded (binary) solutions, respectively. Solid and dashed curves illustrate the trade-off with snapshots of the resulting binary structures for $ka=0.5$. At $ka=1$ (dashed lines), resonant designs are readily achieved even at high~$\Sf$ values, with thresholding having minimal impact on performance. However, at $ka=0.5$, significant structural modifications ($\Sf \leq 0.8$) are necessary to approach resonance, consistent with the physical limits of electrically small antennas~\cite{Capek_etal_2019_OptimalPlanarElectricDipoleAntennas}. For electrically smaller regions, thresholding significantly affects performance as sharp resonances are sensitive to even slight geometric modifications.

It's worth noting that despite similar performance in $\lambda_1$, diverse solutions are possible within the resonance constraint. The surface area constraint provides a valuable metric for comparing the performance of different resonant structures, allowing designers to optimize performance within given size restrictions.

\subsection{Resonance Tuning With Improved Modal Bandwidth}
We obtained different resonant structures under various surface area constraints in the previous section. While these solutions successfully achieve resonance of the first characteristic mode, each solution inherently possesses different modal fractional bandwidth potential, which is usually quantified by modal radiation Q-factor~$Q_n$~\cite{HarringtonMautz_ControlOfRadarScatteringByReactiveLoading} and can be used as another selector. The modal Q-factor is defined in the vicinity of the resonance as 
\begin{equation}
    Q_{n} = \dfrac{\M{I}_n^\herm \M{X}_0' \M{I}_n}{2\M{I}_n^\herm \M{R}_0 \M{I}_n} \approx \dfrac{\omega}{2} \dfrac{\partial \lambda_n}{\partial \omega},
    \label{eq:modalQ}
\end{equation}
where matrix~$\M{X}_0'$ is the stored energy matrix~\cite{Vandenbosch_ReactiveEnergiesImpedanceAndQFactorOfRadiatingStructures} defined as
\begin{equation}
\M{X}_0' = \omega\dfrac{\partial \M{X}_0}{\partial \omega}.
\label{eq:stored matrix}
\end{equation}
The modal Q-factor~\eqref{eq:modalQ} can be included in the optimization as
\begin{equation}
    \begin{aligned}
        \underset{\des}{\T{minimize}} &\quad Q_{1} + \gamma\left(|\lambda_1|^2 +\nu |\delta_1|^2\right) \\
         \T{subject~to} 
         &\quad \left(\M{X}_0 - \J \M{R}_\rho(\desF)\right)\M{I}_1 = \xi_1 \M{R}_0 \M{I}_1,\\
          &\quad \desF = \heavisideFilter({\densityFilter({\des}})),\\
         &\quad 0\leq \des \leq 1,\quad t=1,\dotsc,T,\\
    \end{aligned}
    \label{eq:minQ}
\end{equation}
where~$\gamma$ is a penalization parameter that forces the design into the resonance, and~$\nu$ is the parameter introduced in~\secref{sec:resonance}. As stated, the goal is to make the first significant characteristic mode resonant while optimizing its characteristic current to achieve the minimal value of the modal Q-factor. Another option would be to combine more modes to reach the resonance~\cite{CapekJelinek_OptimalCompositionOfModalCurrentsQ}, but this is out of the scope of this paper.

\begin{figure}[!t]
\centering
\includegraphics[width=3.25in,clip]{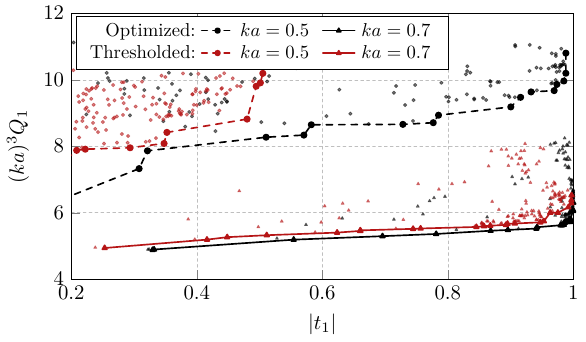}
\caption{Multiple topology optimization runs~\eqref{eq:minQ} showing achieved modal significance~$|t_1|$ and associated modal Q-factor~$Q_1$ of the first characteristic mode. Results are compared at two electrical sizes $ka = 0.5$ and $ka = 0.7$, where black markers represent optimized gray designs and red markers show their thresholded binary counterparts. Non-dominated solutions are connected with dashed and solid lines, respectively. Higher values of parameter~$\gamma$ in~\eqref{eq:minQ} enforce resonance during the optimization ($|t_1|\to1$). Thresholding will inevitably shift the resonance from the target electrical size.}
\label{fig:Quality factors}
\end{figure}

The optimization is performed on a rectangular design region as in the previous section for electrical sizes $ka=\{0.5,0.7\}$ using the RCI cluster~\cite{rci}. Multiple runs explore different parameters~$\gamma\in\left[0.1,100\right]$ and~$\nu=\{0.1,0.5\}$. Each iteration requires 96~seconds on average, a significantly longer time than optimizing characteristic numbers in~\secref{sec:resonance}. This is due to the calculations of the adjoint system of size $(2N+2) \times (2N+2)$, see~\appref{app:adjoint}.

The individual points show the outcome of each topology optimization run in~\figref{fig:Quality factors} for both electrical sizes. The black and red markers represent optimized (gray) and thresholded (binary) solutions.
The modal significance~$|t_n|$~\cite{cabedo2007theory} of a characteristic mode, defined as
\begin{equation}
|t_n| = \dfrac{1}{ 1 + \lambda_n^2 },
\label{eq:modal_significance}
\end{equation}
is utilized to quantify how well a mode can be excited, \ie{}, markers in~\figref{fig:Quality factors} show excitability of a mode and its potential fractional bandwidth through the associated modal Q-factor which is normalized with~$(ka)^3$, as this scaling factor relates to the bandwidth limitations of electrically small antennas~\cite{JonssonGustafsson_StoredEnergiesInElectricAndMagneticCurrentDensities_RoyA}. Due to the scalar nature of the objective and the framework's local nature, not every optimization run yields a non-dominated solution, suggesting the need to rerun the optimization multiple times. The dashed and solid lines illustrate the trade-off between modal Q-factor and modal significance for both electrical sizes, connecting non-dominated solutions obtained through parameter settings in~\eqref{eq:minQ}. Although resonant solutions ($|t_1|\to1$) are achievable for high values of~$\gamma$, the subsequent thresholding procedure shifts the resonance and unpredictably changes the modal Q-factor. This shift is evident from the separation between black markers (optimized designs) and red markers (thresholded designs) in~\figref{fig:Quality factors}. Furthermore, thresholding impacts performance more significantly in electrically smaller regions, which are more sensitive to structural changes.
 
Two representative thresholded designs are selected for detailed analysis in~\tabref{tab:Qfactor}, where they are compared with an empirical meander line antenna design~\cite{Capek_etal_2019_OptimalPlanarElectricDipoleAntennas} which serves as a benchmark due to its common applications~\cite{Best_ElectricallySmallResonantPlanarAntennas}. The lower bound~$Q_\T{lb}^{\T{TM}}$~\cite{GustafssonTayliEhrenborgEtAl_AntennaCurrentOptimizationUsingMatlabAndCVX} restricted to the radiation of the sole transverse magnetic mode is also included to assess the optimality of the designs. At $ka=0.5$, the optimization cannot match the meanderline performance shown in~\tabref{tab:Qfactor}. While a finer grid with a smaller density filter radius could theoretically yield better-performing solutions, this approach was not pursued due to computational limitations. However, the topology optimization achieves superior performance at $ka=0.7$, where meander line geometry simplifies. 

\begin{table}[!t] 
\centering
\caption{Modal significance~$|t_1|$ and associated modal Q-factor~$Q_1$ of the first characteristic mode of thresholded designs from~\figref{fig:Quality factors} and of empirical meanderline antenna~\cite{Capek_etal_2019_OptimalPlanarElectricDipoleAntennas}. The optimality of designs is assessed on the distance to the lower bound~$Q_\T{lb}^{\T{TM}}$~\cite{GustafssonTayliEhrenborgEtAl_AntennaCurrentOptimizationUsingMatlabAndCVX}.}
\renewcommand{\tabcolsep}{0.20cm}

\newcolumntype{E}{ >{\centering\arraybackslash} m{0.5cm} }
\newcolumntype{A}{ >{\centering\arraybackslash} m{2.5cm} }
\newcolumntype{B}{ >{\centering\arraybackslash} m{1.25cm} }
\newcolumntype{C}{ >{\centering\arraybackslash} m{1.25cm} }
\newcolumntype{D}{ >{\centering\arraybackslash} m{1.25cm} }

\begin{tabular}{ E A B C D} 
& Structure & $|t_1|$ & $(ka)^3Q_1$ & $(ka)^3Q_\T{lb}^\T{TM}$\\[1pt] \toprule
\multirow{5}{*}{{\rotatebox[origin=c]{90}{$ka=0.5$}}}&\includegraphics[width=2.5cm]{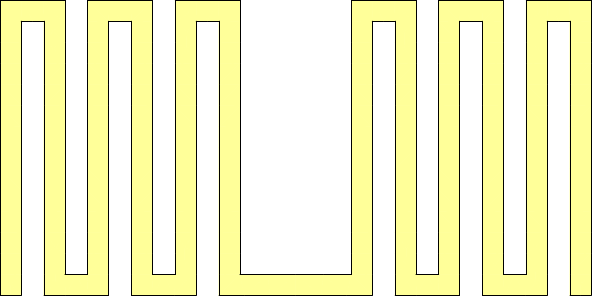} & 1.00 & 5.83 & \multirow{5}{*}{5.27} \\
&\includegraphics[width=2.5cm]{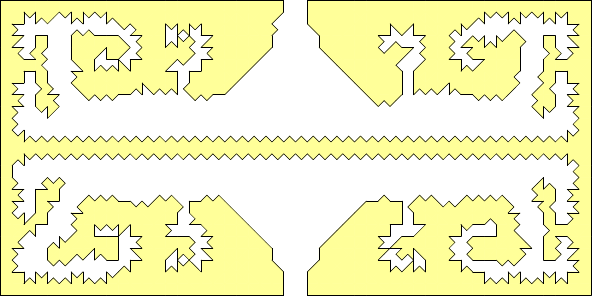} & 0.49 & 8.83 &  \\\midrule
\multirow{5}{*}{{\rotatebox[origin=c]{90}{$ka=0.7$}}}&\includegraphics[width=2.5cm]{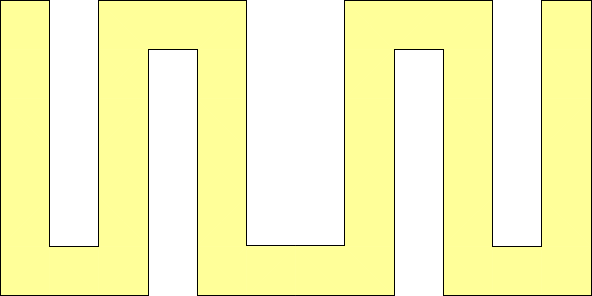} & 1.00 & 6.12 & \multirow{5}{*}{5.25} \\
&\includegraphics[width=2.5cm]{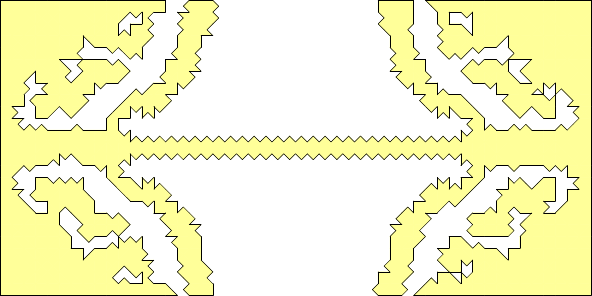} & 0.91 & 5.74 &  \\
\bottomrule
\end{tabular}
\label{tab:Qfactor}
\end{table}

\subsection{Simultaneous Manipulation with Several Modes}
The orthogonality of characteristic mode radiation patterns provides an effective foundation for MIMO antenna design~\cite{li2014designHandsetAntennas},~\cite{dong2020design}. The key goal in MIMO design using CMA is to obtain several resonant characteristic modes simultaneously. Such multi-mode resonance task can be formulated as a \Quot{min-max} problem
\begin{equation}
    \begin{aligned}
        \underset{\des}{\T{minimize}}  &\quad \T{max}\left\{|\lambda_1|,\dots,|\lambda_n| \right\},  
        \label{eq:minmax}
    \end{aligned}
\end{equation}
which is further written as bound formulation~\cite{BendsoeSigmund_TopologyOptimization}
\begin{equation}
    \begin{aligned}
        \underset{\des}{\T{minimize}}  &\quad z \\
         \T{subject~to} 
         &\quad |\lambda_n|^2 +\nu |\delta_n|^2 \leq z\\
         &\quad \left(\M{X}_0 - \J \M{R}_\rho(\desF)\right)\M{I}_n = \xi_n \M{R}_0 \M{I}_n,\\
          &\quad \desF = \heavisideFilter({\densityFilter({\des}})),\\
         &\quad 0\leq \des \leq 1,\quad t=1,\dotsc,T,\\
         &\quad z\geq0,\\
         &\quad n\in\{1,2,3\},
         \label{eq:CMA MIMO}
    \end{aligned}
\end{equation}
where an auxiliary variable~$z$ is introduced to handle the non-differentiable max function. The optimization task~\eqref{eq:CMA MIMO} can be efficiently solved using MMA~\cite{svanberg2007mma}. 

\begin{figure}[!t]
\centering
\includegraphics[width=3.25in]{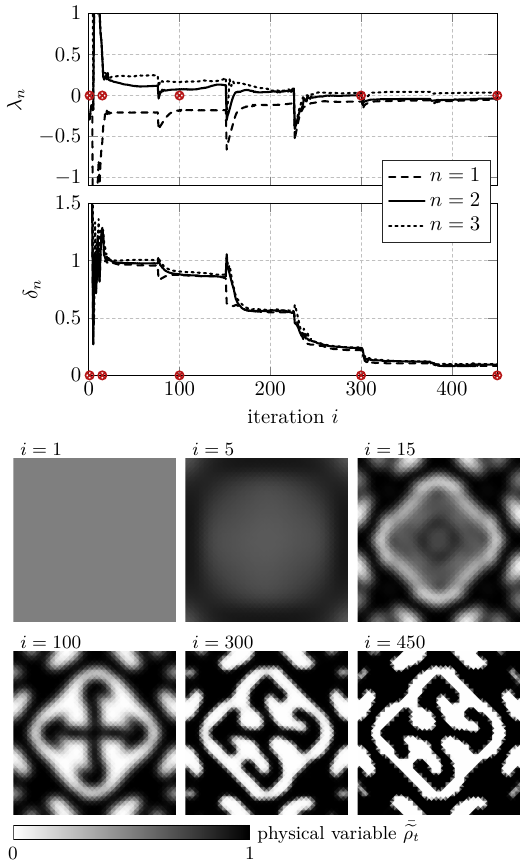}
\caption{The convergence history of the topology optimization task~\eqref{eq:CMA MIMO} for three most significant characteristic modes~$n=\{1,2,3\}$. A few snapshots of the density distribution from selected iterations of the run are included below the convergence plots to demonstrate the progression from homogeneous distribution ($i=1$) to the final optimized structure ($i=450$).}
\label{fig:Convergence_MIMO}
\end{figure}

A critical aspect of this formulation lies in the adjoint sensitivity analysis, see~\appref{app:adjoint}, which assumes distinct characteristic numbers only. This assumption becomes theoretically violated at the optimal solution where all modes achieve resonance ($\lambda_n=0$ for each~$n$). However, this theoretical limitation has a minimal practical impact on the proposed optimization process since the characteristic numbers remain distinct throughout most of the optimization due to intermediate densities distributed in a non-symmetric way, and even if degeneracy occurs, it affects the gradient computation only temporarily. 

The optimization task~\eqref{eq:CMA MIMO} is solved at electrical size~$ka=2$ with $\nu=0.5$ for the three most significant modes ($n={1,2,3}$) on a square region uniformly discretized into 4900 triangles (7280 basis functions). Each iteration requires approximately 110 seconds.

The convergence plot shown in~\figref{fig:Convergence_MIMO} demonstrates the evolution of both parts of the optimized characteristic numbers~$\xi_n$. The convergence curve is stable despite the initial degeneracy between the first two modes, which results in a slight oscillation. However, modes maintain their distinctness while approaching the resonance.

\begin{figure}[!t]
\centering
\includegraphics[width=3.25in]{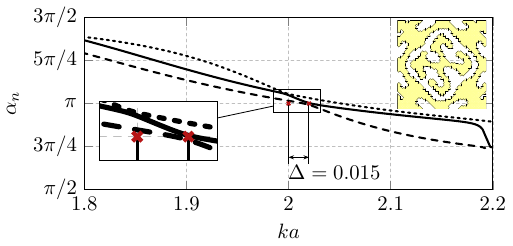}
\caption{Characteristic angles~$\alpha_n$ versus electrical size $ka$ for three characteristic modes of the thresholded binary design. The resonance point is shifted from the target electrical size~$ka=2$. The non-crossing behavior of modes follows the von Neumann-Wigner theorem~\cite{vonNeumannWigner_OnTheBehaviourOfEigenvaluesENG} due to the structure's non-symmetry~\cite{SchabEtAl_EigenvalueCrossingAvoidanceInCM}.}
\label{fig:thresholding_MIMO}
\end{figure}

The optimization converges to a structure combining two crossed dipoles and a loop, with all three modes close to the resonance. The snapshots of the density distribution in~\figref{fig:Convergence_MIMO} show the evolution from homogeneous distribution to the final optimized structure. After thresholding of the optimized solution~\cite{tucek2023TopOptMoM_minQ}, the resonances of the modes are slightly shifted from the target value $ka$, see the characteristic angles in~\figref{fig:thresholding_MIMO}. The non-symmetric binary structure exhibits mode behavior consistent with the von Neumann-Wigner theorem~\cite{vonNeumannWigner_OnTheBehaviourOfEigenvaluesENG}.

\begin{figure}[!t]
\centering
\includegraphics[width=3.25in]{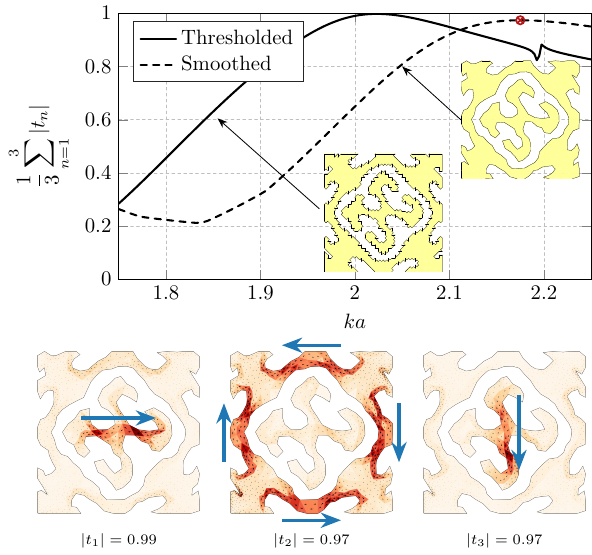}
\caption{The sum of the modal significances~$t_n$ of the three characteristic modes of the thresholded binary design and the design with smoothed boundaries. The characteristic current densities of the three most significant modes are included below for the smoothed structure at resonant electrical size~$ka=2.17$. Blue arrows indicate the spatial current flow.}
\label{fig:ModalSignificance_MIMO}
\end{figure}

The hard-thresholding results in a binary design with jagged boundaries due to underlying triangular discretization, see the design in~\figref{fig:thresholding_MIMO}. The shape optimization techniques~\cite{BendsoeSigmund_TopologyOptimization} are commonly used within FEM~\cite{wang2017antenna} to refine such boundaries, but these approaches have not yet been established for MoM models. Therefore, we address this issue by applying a Gaussian convolution filter~\cite{gonzalez2009digital} for blurring the structure's boundary to eliminate complicated features. As shown in~\figref{fig:ModalSignificance_MIMO}, the smoothing process maintains high levels of modal significance for all modes but results in a minor performance drop and a resonance shift from the target electrical size to~$ka = 2.17$. However, the resonance change can be compensated by rescaling the structure.
The characteristic current densities, see~\figref{fig:ModalSignificance_MIMO}, illustrate the spatial distribution of the three most significant modes of the smoothed structure at its resonant electrical size $ka = 2.17$. The spatial separation of the current maxima for each mode suggests their potential excitability~\cite{cabedo2007theory}, though practical implementation may require additional fine-tuning of feeding positions. 

The proposed method can handle multiple-mode manipulation despite the inherent limitations of the adjoint sensitivity analysis. The presented results are especially interesting for the design of MIMO antenna structures supporting multiple resonant characteristic modes.

\section{Conclusion}\label{sec:conclusion}
The density-based topology optimization was formulated with characteristic mode analysis to optimize arbitrary modal quantities of conductor-based structures under method-of-moments modeling. Introducing a material matrix as an auxiliary variable facilitated density formulation, albeit making characteristic numbers and currents complex. Adjoint sensitivity analysis was derived for distinct characteristic numbers, requiring the solution of $(N+1)\times(N+1)$ or $(2N+2)\times(2N+2)$ systems depending on the objectives. The proposed method employed a continuation approach with projection and density filters to achieve near-binary designs with an imposed length-scale in the weak sense. The approach demonstrated superior computational efficiency compared to alternative formulations and can handle multi-frequency optimization with only a linear increase in complexity. While the method delivers near-binary solutions, postprocessing via thresholding can lead to performance drops and resonance shifts, which may be compensated by rescaling the structure or combining it with binary topology optimizers. The method effectively separates geometry and feeding synthesis, potentially inspiring novel designs with \textit{a posteriori} feeding synthesis. 

The method's efficiency and properties have been demonstrated in single- and multiple-mode manipulations with characteristic numbers. Due to self-penalization properties, penalizing the imaginary part of the characteristic number can yield better solutions and directly influence the performance of the resulting binary structure. The optimization settings are not universal and must be tailored to a particular design goal. However, the superior convergence of the method allows for multiple runs with various settings and heuristic parameter values to obtain the best possible solution.

Future research possibilities include extending the framework pool to new application problems, \eg{}, selective mode optimization, combining more modes, etc. Addressing the jagged boundaries of solutions through automated smoothing while preserving performance values remains an area for improvement. Extending the adjoint sensitivity analysis to include degenerated characteristic numbers could enhance the utility of the method in antenna design, where degeneracy is often leveraged.


%

\appendices
\section{Influence of the Surface Resistivity of the Air on the First Characteristic Mode}
\label{app:Interpolation behavior Zs}

The influence of vacuum surface resistivity settings~$\Zair$ in the interpolation scheme~\eqref{eq:Interpolation} is examined for evaluating the first characteristic number~$\xi_1=\lambda_1 -\J\delta_1$ of the slot-loaded plate at $ka=0.7$, while the density~$\rho$ is varied within the design region of the structure, as shown in~\figref{fig:structure for material model investigations}. This is similar to the investigation performed in~\secref{sec:Influence}, but the electrical size is kept fixed, and the effect of~$\Zair$ on the material model is presented. The $\Zmet=10^{-2}\;\Omega$ is kept fixed.

\begin{figure}[!t]
\centering
\includegraphics[width=3.25in,clip]{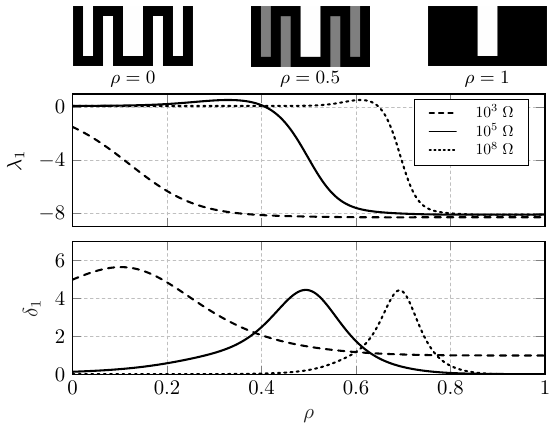}
\caption{The influence of varying surface resistivity~$\Zair$ associated in topology optimization with the vacuum. The interpolation function~\eqref{eq:Interpolation} is examined in the context of evaluating the first characteristic number~$\xi_1=\lambda_1-\J\delta_1$ of meander line antennas at $ka=0.7$, see~\figref{fig:structure for material model investigations}. The structure is parametrized by density function~$\rho$ defined within the design region and representing surface resistivity~$\Rs(\rho)$ in~\eqref{eq:Interpolation}, see snapshots of the structure while $\rho$ is swept. The surface resistivity of the conductor is kept fixed at~$\Zmet=10^{-2}\;\Omega$ while $\Zair$ is varied.}
\label{fig:Interpolation Behavior Zs}
\end{figure}

Both parts of the first characteristic number~$\lambda_1$ and~$\delta_1$ are evaluated while density~$\rho$ in the design region is swept for three settings of~$\Zair$, see~\figref{fig:Interpolation Behavior Zs}. If the density is set to~$\rho=0$, the resonant meander line antenna with $\lambda_1\approx0$ and $\delta_1\to0$ is obtained, see snapshots above the axis in~\figref{fig:Interpolation Behavior Zs}. Furthermore, setting~$\rho=1$ leads to a non-resonant slot-loaded plate. These boundary cases correspond to the investigation performed in~\secref{sec:Influence} since the resonance of the slot-loaded plate is located at $ka\approx1$.    

The dotted and solid lines in~\figref{fig:Interpolation Behavior Zs} represent the physically accurate resistivity of air since they correctly reflect the values of~$\lambda_1$ and~$\delta_1$ for~$\rho=0$ and~$\rho=1$ as opposed to the dashed line which does not. The dotted and solid lines show the most notable differences in the mid-range of~$\rho$. Setting $\Zair=10^5 \;\Omega$ allows the optimizer to explore a broader range of design variable values since a change in the design variable has a more pronounced effect on $\lambda_1$ and $\delta_1$, simplifying an optimization at the cost of a small inaccuracy in the objectives~\cite{tucek2023TopOptMoM_minQ}. Unfortunately, this material model isn't universal, and numerical investigations are always necessary to determine the most suitable form for each specific problem.

\section{Adjoint Sensitivity Analysis for CMA}\label{app:adjoint}
Consider the minimization task of an objective function~$f(\M{I}_n,\xi_n, \rho)$, where the governing equation is CMA formulation~\eqref{eq:lossy CMA}. The adjoint sensitivity analysis~\cite{tortorelli1994design} efficiently computes a gradient information~$\D f/ \D\rho$. However, due to the complex nature of~\eqref{eq:lossy CMA} the characteristic currents are not unique since they can be scaled or rotated in the complex plane~\cite{murthy1988derivatives}. Hence, additional constraints are introduced to ensure the uniqueness of the characteristic currents~\cite{he2023eigenderivation}. The unit radiated power is constrained with bilinear form
\begin{align}
    \M{I}_n^\herm \M{R}_0 \M{I}_n = 1,
    \label{eq:normalization}
\end{align}
and the angle is constrained by forcing one entry of $\M{I}_n$ with the maximum norm to be a real number 
\begin{equation}
    \begin{aligned}
    \T{Im}\left(\M{I}_{n,k}\right) &= 0,\\
    k &= \T{arg} \underset{j}{\max}||\M{I}_{n, j}||_2,
    \label{eq:angle constraint}
\end{aligned}
\end{equation}
where $_j$ is the entry index. The adjoint sensitivity analysis is initiated by adding the residual to the objective function and introducing the augmented Lagrangian as

\begin{align}
    L =& f(\M{I}_\T{Re},\M{I}_\T{Im},\lambda,\delta,\rho) \nonumber\\ &+\M{z}_{\T{Re}}^\trans\left[\left(\M{X}_0 - \lambda\M{R}_0\right)\M{I}_{\T{Re}} + \left(\M{R}_\rho(\rho) - \delta\M{R}_0\right)\M{I}_{\T{Im}}\right] \nonumber\\
    &+\M{z}_{\T{Im}}^\trans\left[\left(-\M{R}_\rho(\rho) + \delta\M{R}_0\right)\M{I}_{\T{Re}} + \left(\M{X}_0 - \lambda\M{R}_0\right)\M{I}_{\T{Im}}\right] \label{eq:Lagrangian}\\
&+a\left[\M{I}_\T{Re}^\trans\M{R}_0\M{I}_\T{Re} + \M{I}_\T{Im}^\trans\M{R}_0\M{I}_\T{Im} - 1\right]\nonumber\\
    &+b \M{e}_k^\trans \M{I}_{\T{Im}}\nonumber,
\end{align}
where $\M{z}_{\T{Re}}$,~$\M{z}_{\T{Im}}$,~$a$,~$b$ are Lagrange multipliers, also called adjoint variables, $\M{e}_k^\trans$ is the unit vector with only non-zero entry at $k$-th position, and for the sake of clarity subscript $_n$ is omitted. Performing a total derivative of~\eqref{eq:Lagrangian} with respect to the design variable~$\rho$ and rearranging the terms yields

\begin{equation}
\begin{aligned}
    \frac{\D L}{\D \rho} &= \frac{\D f}{\D \rho} = \frac{\partial f}{\partial \rho} +\M{z}_{\T{Re}}^\trans\frac{\partial\M{R}_\rho(\rho)}{\partial\rho}\M{I}_\T{Im} - \M{z}_{\T{Im}}^\trans\frac{\partial\M{R}_\rho(\rho)}{\partial\rho}\M{I}_\T{Re}\\
    &+\left[\frac{\partial f}{\partial \lambda} - \M{z}_\T{Re}^\trans \M{R}_0\M{I}_\T{Re} - \M{z}_\T{Im}^\trans \M{R}_0\M{I}_\T{Im}\right]\frac{\partial\lambda}{\partial\rho}\\
    &+\left[-\frac{\partial f}{\partial \delta} + \M{z}_\T{Re}^\trans \M{R}_0\M{I}_\T{Im} - \M{z}_\T{Im}^\trans \M{R}_0\M{I}_\T{Re}\right]\frac{\partial\delta}{\partial\rho}\\
    &+\bigg[\frac{\partial f}{\partial \M{I}_\T{Re}} + \M{z}_\T{Re}^\trans\left(\M{X}_0-\lambda\M{R}_0\right) - \M{z}_\T{Im}^\trans\left(\M{R}_\rho(\rho)-\delta\M{R}_0\right)\\
    &\qquad\qquad+2a\M{I}_\T{Re}^\trans\M{R}_0\bigg]\frac{\partial\M{I}_\T{Re}}{\partial \rho}\\
    &+\bigg[\frac{\partial f}{\partial \M{I}_\T{Im}} + \M{z}_\T{Re}^\trans\left(\M{R}_\rho(\rho)-\delta\M{R}_0\right) + \M{z}_\T{Im}^\trans\left(\M{X}_0-\lambda\M{R}_0\right)\\
    &\qquad\qquad+2a\M{I}_\T{Im}^\trans\M{R}_0 + b \M{e}_k^\trans\bigg]\frac{\partial\M{I}_\T{Im}}{\partial \rho},
    \label{eq:gradient general}
\end{aligned}
\end{equation}
where we conveniently grouped the computationally expensive terms~$\partial \lambda/\partial \rho$, $\partial \delta/\partial \rho$, $\partial \M{I}_\T{Re}/\partial \rho$, $\partial \M{I}_\T{Im}/\partial \rho$. 
These problematic derivatives are eliminated by solving the adjoint equation
\begin{equation}
    \M{A}\M{w} = \M{F},
    \label{eq:Adjoint equation}
\end{equation}
where
\begin{equation}
\M{A} = 
    \begin{bmatrix}
\M{X}_0-\lambda\M{R}_0 & -\M{R}_\rho(\rho)+\delta\M{R}_0 & 2\M{R}_0\M{I}_\T{Re} & \M{0} \\[0.5em]
\M{R}_\rho(\rho)-\delta\M{R}_0 & \M{X}_0-\lambda\M{R}_0 & 2\M{R}_0\M{I}_\T{Im} & \M{e}_k \\[0.5em]
-\M{I}_\T{Re}^\trans\M{R}_0 & -\M{I}_\T{Im}^\trans\M{R}_0 & 0 & 0 \\[0.5em]
\M{I}_\T{Im}^\trans\M{R}_0 & -\M{I}_\T{Re}^\trans\M{R}_0 & 0 & 0
\end{bmatrix}
\end{equation}
and
\begin{equation}
    \begin{aligned}
 \M{w} &= \begin{bmatrix} \M{z}_\T{Re}^\trans & \M{z}_\T{Im}^\trans & a & b \end{bmatrix}^\trans, \\
 \M{F} &= -\begin{bmatrix} \dfrac{\partial f}{\partial \M{I}_\T{Re}} & \dfrac{\partial f}{\partial \M{I}_\T{Im}} & \dfrac{\partial f}{\partial \lambda} & -\dfrac{\partial f}{\partial \delta} \end{bmatrix}^\trans. 
\end{aligned}
\end{equation}

If the~\eqref{eq:Adjoint equation} is solved for $\M{w}$, the gradient expression~\eqref{eq:gradient general} is reduced to
\begin{equation}
    \frac{\D f}{\D \rho} = \frac{\partial f}{\partial \rho} +\M{z}_{\T{Re}}^\trans\frac{\partial\M{R}_\rho(\rho)}{\partial\rho}\M{I}_\T{Im} - \M{z}_{\T{Im}}^\trans\frac{\partial\M{R}_\rho(\rho)}{\partial\rho}\M{I}_\T{Re}.
\end{equation}
Hence, the adjoint equation of size $(2N+2)\times(2N+2)$ needs to be evaluated once for each iteration of the topology optimization and for each extra constraint to compute the corresponding sensitivities. This derivation holds only for distinct characteristic numbers. If multiple characteristic numbers coincide, the matrix $\M{A}$ in \eqref{eq:Adjoint equation} becomes singular.

Furthermore, if the characteristic current is not involved in the objective function evaluation, the extra constraint~\eqref{eq:angle constraint} is not necessary~\cite{lee2007adjoint} and the adjoint equation~\eqref{eq:Adjoint equation} is reduced to
\begin{equation} 
    \begin{bmatrix}
        \M{X}_0 -\J\M{R}_\rho(\rho) - \xi\M{R}_0 & -\M{R}_0\M{I}\\
        -\M{I}^\herm\M{R}_0 & 0
    \end{bmatrix}
    \begin{bmatrix}
        \M{z}\\
        a
    \end{bmatrix}
    =
    \begin{bmatrix}
        \M{0}\\
        -\dfrac{\partial f}{\partial \xi}
    \end{bmatrix},
    \label{eq:adjoint lambda}
\end{equation}
where $\M{I} = \M{I}_\T{Re} + \J \M{I}_\T{Im}$, and $\M{z} = \M{z}_\T{Re} + \J \M{z}_\T{Im}$. Hence, the adjoint equation is reduced to the size $(N+1)\times (N+1)$ but is composed of complex values. 

If an objective or constraint is dependent both on~$\xi_n$ and~$\M{I}_n$, the adjoint equation~\eqref{eq:Adjoint equation} of size $(2N+2)\times (2N+2)$ must be solved every iteration to obtain gradient information. However, if only~$\xi$ is involved in the evaluation of an objective, only the equation~\eqref{eq:adjoint lambda} of size $(N+1)\times (N+1)$ needs to be solved. 



\ifCLASSOPTIONcaptionsoff
  \newpage
\fi

\begin{IEEEbiography}
[{\includegraphics[width=1in,height=1.25in,clip,keepaspectratio]{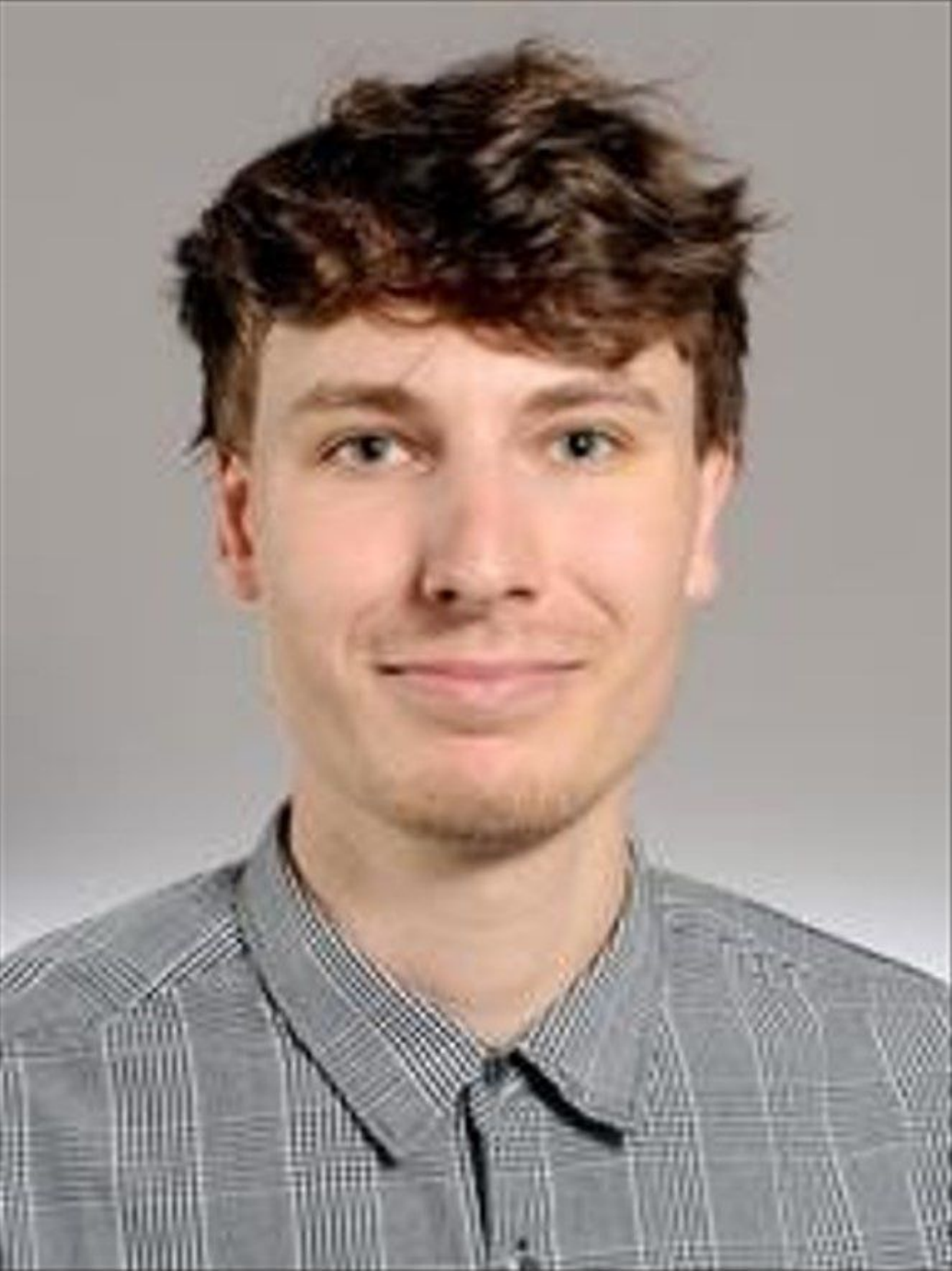}}]{Jonas Tucek}
received an M.Sc. degree in electrical engineering from Czech Technical University, Prague, Czech Republic, in 2021, where he is currently pursuing a Ph.D. degree in the area of topology optimization in electromagnetics. 
\end{IEEEbiography}

\begin{IEEEbiography}[{\includegraphics[width=1in,height=1.25in,clip,keepaspectratio]{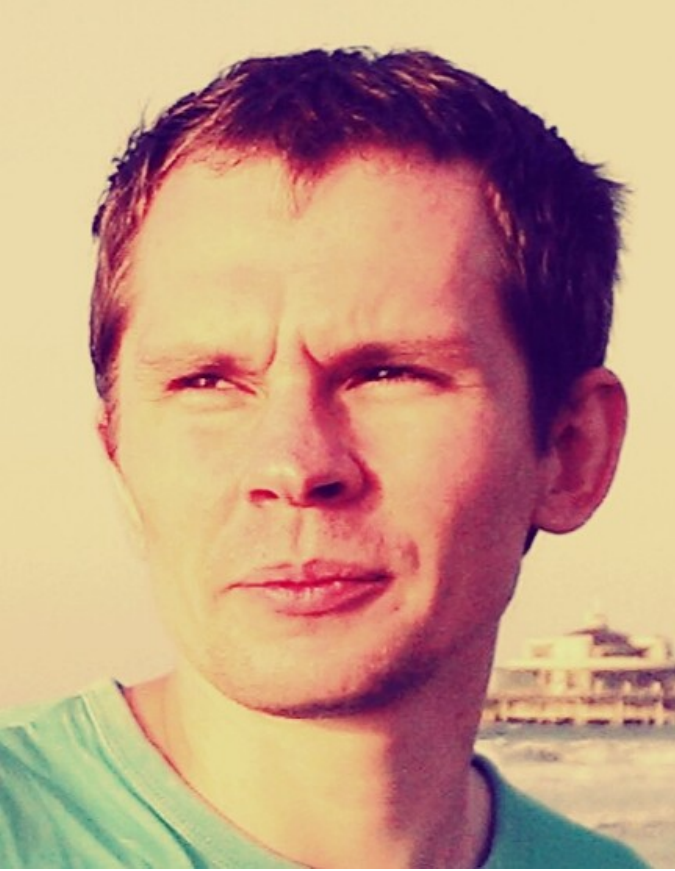}}]{Miloslav Capek}
(M'14, SM'17) received the M.Sc. degree in Electrical Engineering 2009, the Ph.D. degree in 2014, and was appointed a Full Professor in 2023, all from the Czech Technical University in Prague, Czech Republic.
	
He leads the development of the AToM (Antenna Toolbox for Matlab) package. His research interests include electromagnetic theory, electrically small antennas, antenna design, numerical techniques, and optimization. He authored or co-authored over 170~journal and conference papers.

Dr. Capek is the Associate Editor of IET Microwaves, Antennas \& Propagation. He received the IEEE Antennas and Propagation Edward E. Altshuler Prize Paper Award~2022 and ESoA (European School of Antennas) Best Teacher Award in~2023.
\end{IEEEbiography}

\begin{IEEEbiography}[{\includegraphics[width=1in,height=1.25in,clip,keepaspectratio]{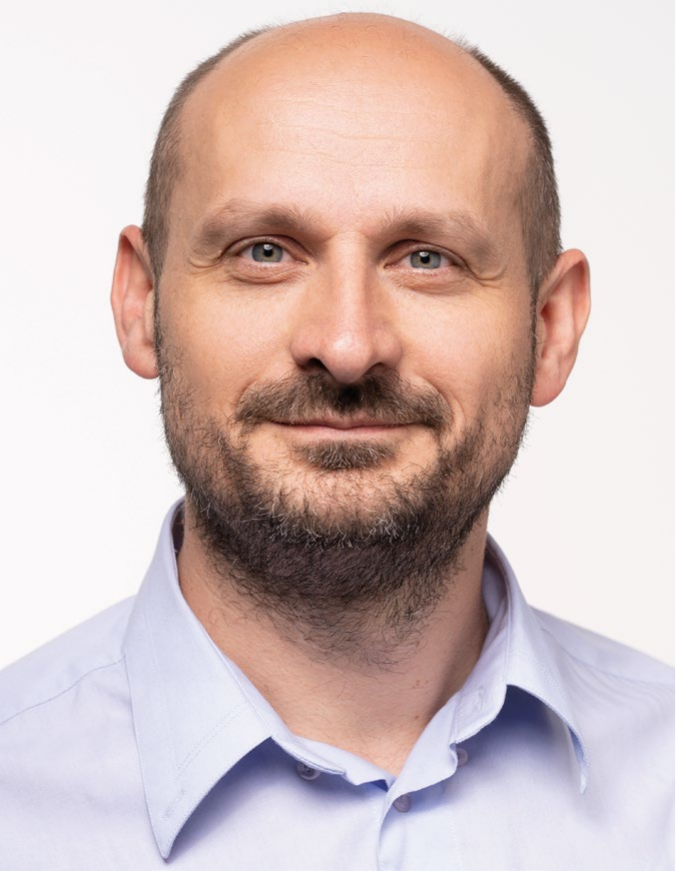}}]{Lukas Jelinek} was born in the Czech Republic in 1980. He received his Ph.D. from the Czech Technical University in Prague, Czech Republic, in 2006 for his work in metamaterials. In 2015, he received a permanent position at the Department of Electromagnetic Field at the same university.

His research interests include wave propagation in complex media, electromagnetic field theory, metamaterials, numerical techniques, and optimization.
\end{IEEEbiography}





\end{document}